\documentclass[11pt]{amsart}
\usepackage{standard}
\usepackage{graphicx}
\usepackage{amsrefs}
\usepackage{amssymb}

\DeclareMathOperator{\singsupp}{sing-supp}
\DeclareMathOperator{\LSpec}{L-Spec}
\DeclareMathOperator{\DLSpec}{Diff-L-Spec}
\newcommand{\D}{\mathcal{D}}
\newcommand{\W}{\mathcal{W}}

\newcommand{\restrictedto}{\rvert}
\DeclareMathOperator{\length}{length}
\newcommand{\dom}{\mathcal{D}}
\DeclareMathOperator{\Res}{Res}
\DeclareMathOperator{\Tr}{Tr}
\newcommand{\bJ}{\mathbf{J}}
\let\det=\relax
\DeclareMathOperator{\det}{det}
\newcommand{\wavegen}{\mathcal{A}}

\title[Diffractive Propagation on Conic Manifolds]{Diffractive Propagation on Conic Manifolds\\S\'eminaire Laurent Schwartz\\January 26, 2016}

\thanks{The author is grateful to Richard Melrose for helpful comments
  on the manuscript.  This work was partly supported by NSF grant DMS--1265568.}

\author{Jared Wunsch}

\address{Department of Mathematics\\Northwestern University\\2033
  Sheridan Rd.\\Evanston IL 60208\\USA}

\email{jwunsch@math.northwestern.edu}

\date{\today}

\begin{document}

\begin{abstract}
  In this survey, we review some applications and extensions of the
  author's results with Richard Melrose on propagation of
  singularities for solutions to the wave equation on manifolds with
  conical singularities.  These results mainly concern: the local
  decay of energy on noncompact manifolds with diffractive trapped
  orbits (joint work with Dean Baskin); singularities of the wave
  trace created by diffractive closed geodesics (joint work with
  G.~Austin Ford); and the distribution of scattering resonances
  associated to such closed geodesics (joint work with Luc Hillairet).
\end{abstract}

\maketitle

\section{Introduction}

Consider the wave equation on a Riemannian manifold $X:$
$$
\Box u=0\text{ on } \RR\times X
$$
where $\Box=D_t^2-\Lap_g,$
$$
\Lap_g=\sum \frac{1}{\sqrt{g}} D_j g^{jk} \sqrt{g} D_k
$$
and $D_j\equiv i^{-1} \pa_{x_j}$.

If $X$ happens to be an odd dimensional Euclidean space, then \emph{Huygens' Principle} applies, i.e.,
the solution
$$
\cos t\sqrt{\Lap} \delta_q
$$
which has initial data a delta-function (and initial derivative zero)
is supported exactly on sphere of radius $\abs{t}.$ In even space
dimensions, or on a general odd dimensional manifold, this principle
is well known to fail, but quite a nice proxy for it persists: we in
general have
$$
\singsupp u(t) \subset \big\{p:\ \text{there exists a geodesic of length } \abs{t} \text{ with endpoints }p,q\big\}.
$$
(Recall that the singular support of a distribution is the set of points near which is it not locally a smooth function.)
A more precise result yet is the refinement of this statement to deal
with the \emph{wavefront set} of the distribution $u;$ $\WF u$ is a
conic closed subset of $T^*X$ such that $\pi \WF u=\singsupp u.$   H\"ormander's
rather general theorem \cite{Hormander9} on propagation of
singularities tells us in this special case that for a solution $u$ of
the wave equation, $\WF u$ is invariant
under the (forward and backward) geodesic flow on $T^*X.$  Thus the
initial wavefront given by (the lift to the light cone of)
$N^* \{q\}$ then spreads into the conormal bundle of expanding
distance spheres.

Generalizing this result to manifolds with boundary (with Dirichlet or
Neumann boundary conditions) turns out be a rather complicated story.
Chazarain \cite{Ch:73} showed that singularities striking the boundary
transversely simply reflect according to the usual law of geometric
optics (conservation of energy and tangential momentum, hence ``angle
of incidence equals angle of reflection'') for the reflection of
bicharacteristics.  The difficulties arise, however, in the treatment
of geodesics tangent to the boundary: in \cite{Melrose-Sjostrand1} and
\cite{Melrose-Sjostrand2} Melrose--Sj\"ostrand showed that, at
these ``glancing points,'' singularities may only propagate along certain
generalized bicharacteristics.  By parametrix constructions of Melrose
\cite{Melrose14} and Taylor \cite{Taylor1}, these $\CI$ singularities
do \emph{not} propagate along concave boundaries (e.g.\ they do not
``stick'' to the exterior of a convex obstacle).  Note that this last
result ceases to be true in the analytic, rather than smooth,
category.

\begin{figure}[bth]\label{figure:fundamental}
\includegraphics[scale=1.3,angle=-90]{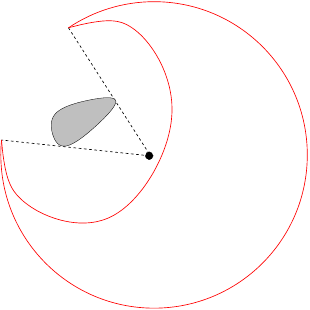}\caption{Singularities
  of the fundamental solution of the wave equation exterior to a convex obstacle.}\end{figure}

A simple summary of some of the fundamental results in the subject is provided by Figure~\ref{figure:fundamental}.
This figure shows the singularities of the fundamental solution the
wave equation in the exterior of a convex obstacle in the plane.
There is (part of) a circular front of directly propagated
singularities as well as a curved front of singularities reflected off
the obstacle in accordance with Snell's law.  Most crucially,
\emph{there are no singularities behind the obstacle} in the ``shadow
region,'' as a consequence of the parametrix construction of Melrose and Taylor.

\begin{figure}
\includegraphics[scale=0.7]{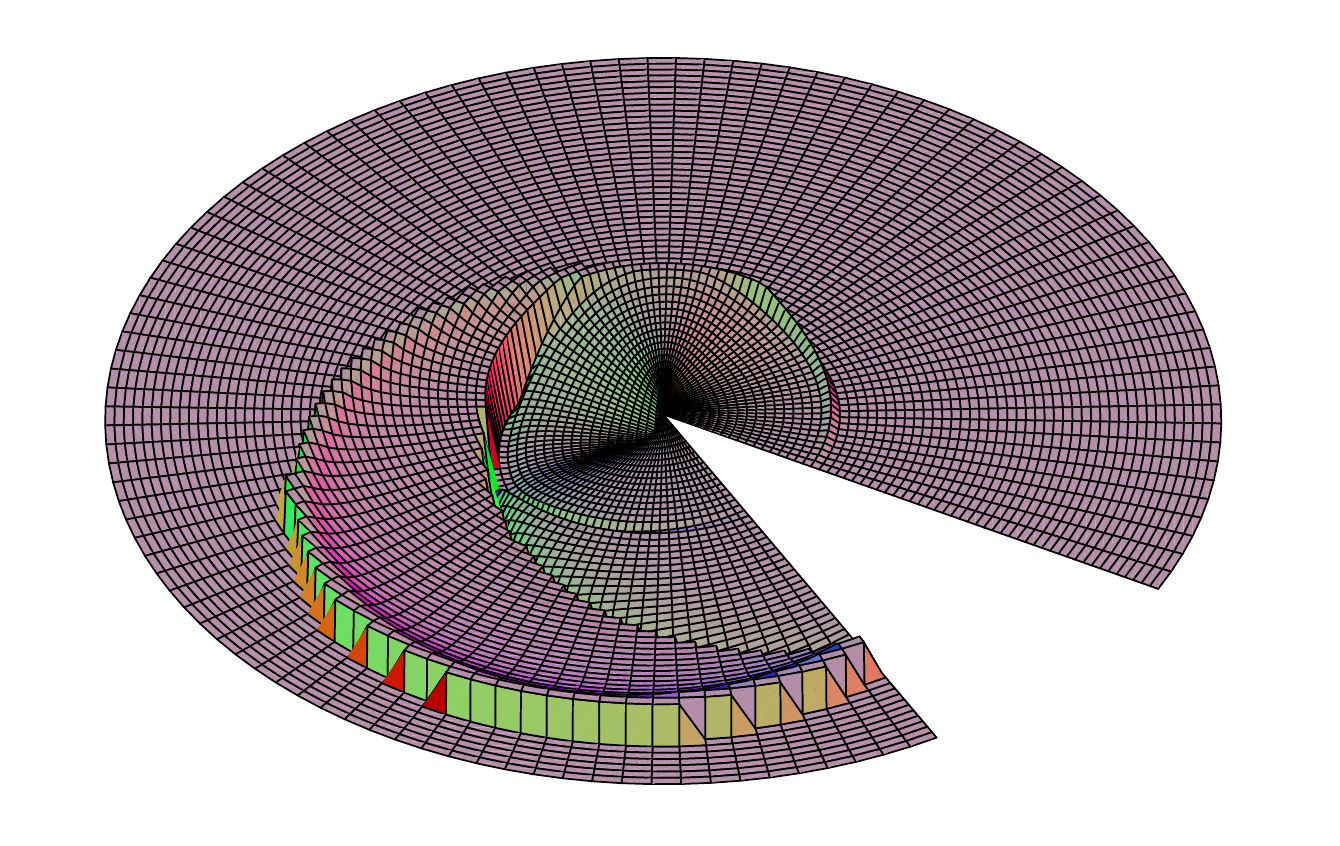}\caption{Singularities
  of the fundamental solution of the wave equation exterior to a
  wedge. \label{figure:wedge}}\end{figure}

By contrast, it has been known since the late 19th century (starting
with work of Sommerfeld \cite{Sommerfeld1}) that if the obstacle has a
sharp corner, singularities \emph{do} propagate, i.e., \emph{diffract,}
into the shadow region behind the obstacle.  Figure~\ref{figure:wedge}
shows the fundemental solution of the wave equation in the exterior of
a wedge; we can easily see a circular wave of singularities emanating
from the tip of the wedge and giving rise to singularities in the
shadow region.

As alluded to above, general boundaries present special difficulties
of their own, so in order to study the diffraction phenomenon in a
simple setting, we now mostly set aside this class of manifolds, and focus on
\emph{manifolds with conic singularities} where wave equation
solutions will exhibit diffraction, but the geometry of geodesics is
relatively manageable.

\section{Conic geometry}

We define a \emph{conic manifold}
to be a manifold $X$ (of dimension $n$) with boundary $Y=\pa X,$ and a Riemannian metric
on $X^\circ$ such that in terms of some boundary defining function $x$
we have in a collar neighborhood of $Y,$
$$
g=dx^2 +x^2 h
$$
where $h$ is a smooth symmetric 2-cotensor such that $h|_{Y}$ is a
metric on $Y.$ Note in particular that $g$ degenerates at $\pa X$ so
as not to be a metric uniformly up to the boundary.  

\begin{figure}[bth]
\includegraphics[scale=.2,angle=-90]{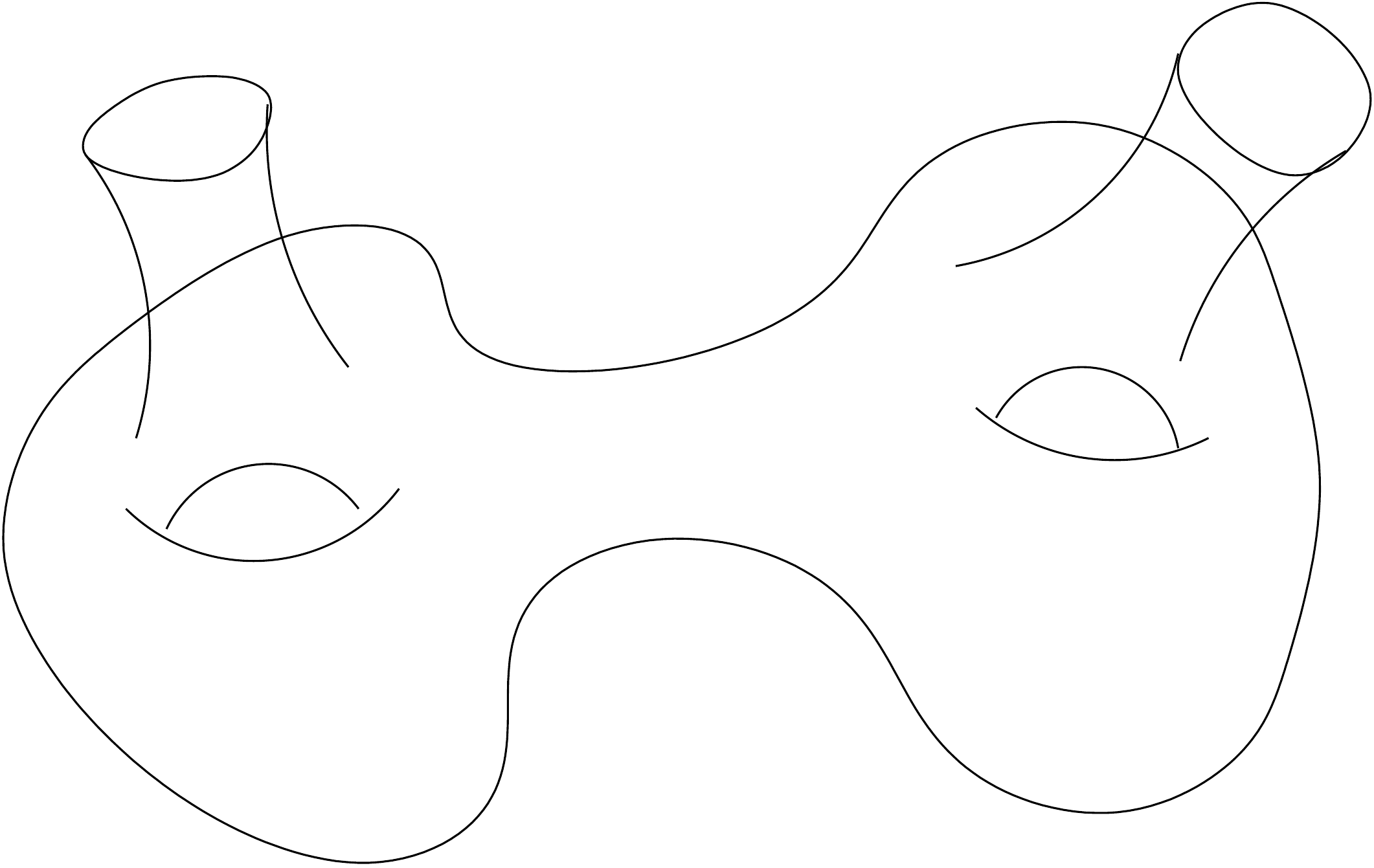}
\includegraphics[scale=.2,angle=-90]{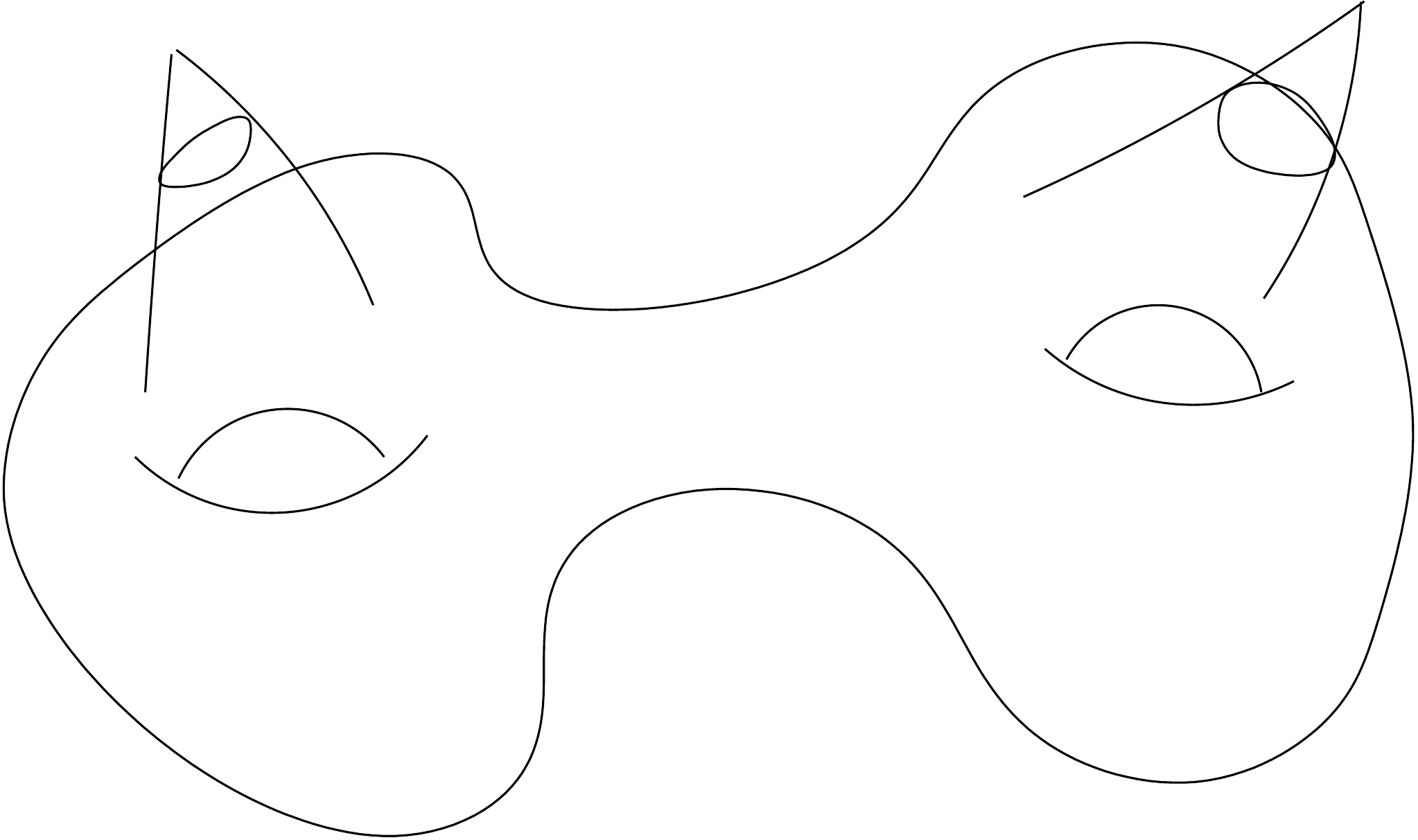}
\caption{Smooth structure, and Riemannian picture of $X$ \label{figure:conicgeometry}}
\end{figure}

The upshot is that while $X$ looks like a manifold with boundary from
the point of view of $\CI$ structure, it is metrically a manifold with
conic singularities: from the point of view of metric geometry, if we
write the connected components of the boundary as 
$$
Y=\bigsqcup Y_i
$$
then each boundary component $Y_i$ should be viewed as a \emph{cone
  point}.  (See Figure~\ref{figure:conicgeometry}.)

The conic manifold as defined here should thus be viewed as a manifold
with conic singularities \emph{already equipped} with the blow-up that
has desingularized it to a smooth manifold with boundary.  Here the
cost of having a smooth manifold is of course having a degenerate metric.

A very special case of a conic manifold is that of a surface
obtained by gluing together two copies of the interior (or exterior)
of a polygonal planar domain along their common edges.  This gives a
flat surface with cone points where the polygon had vertices.  The
study of the wave equation on the original domains with Dirichlet/Neumann
conditions is equivalent to the study of odd/even solutions of the wave equation on the
doubled manifold---see Hillairet \cite{Hillairet:2005}.

The behavior of geodesics on conic manifolds is of considerable
interest near the cone point.  The crucial observation is that it is
in fact quite hard to aim a geodesic so as to hit the cone point: most
will pass nearby and miss.  Indeed, starting out near the cone point,
there is a unique direction to aim in, in order to reach a nearby cone point.
\begin{proposition}[Melrose--Wunsch \cite{Melrose-Wunsch1}]
Every $y \in Y=\pa X$ is the endpoint of a unique geodesic; these
geodesics foliate a collar neighborhood of $Y$:
\end{proposition}
This is equivalent to a normal-form statement for the metric: we can
find coordinates so that $h=h(x,y,dy)$ has no $dx$ components, and
thus the curves $x=x_0 \pm t, y=y_0$ are unit-speed geodesics.

A crucial point in trying to make sense of propagation of singularites
is to make a reasonable definition of the continuation of a geodesic
that reaches a cone point.  There are two reasonable candidates for
this definition, one more restrictive than the other, and both play a
role here:

\begin{definition} We define geodesics passing through $\bigsqcup Y_j \equiv \pa X$ as follows:
\begin{itemize}
\item A \emph{diffractive geodesic} is a geodesic which, upon reaching
  the boundary component $Y_i$ along a geodesic ending at a point
  $y\in Y_i,$ immediately then leaves the boundary from some point $y'
  \in Y_i.$
\item A \emph{geometric geodesic} is a geodesic which, upon reaching
  the boundary component $Y_i$ along a geodesic ending at a point
  $y\in Y_i,$ immediately then leaves the boundary from some point $y'
  \in Y_i$ such that $y,y'$ are endpoints of a geodesic \emph{in $Y_i$}
  (w.r.t.\ the metric $h|_{Y_i})$ of length $\pi.$
\item A \emph{strictly diffractive} geodesic is one which is
  diffractive but not geometric.
\end{itemize}
\end{definition}

A more intuitive definition of geometric geodesics is as follows: they
are the geodesics that are \emph{locally approximable} by families of geodesics in
$X^\circ$.  We refer the reader to \cite{Melrose-Wunsch1} for more
detail on these definitions.

\section{Propagation of singularities on conic manifolds}

Consider now solutions to the wave equation on a manifold with conic
singularities.  We always employ the \emph{Friedrichs extension} of
the Laplacian acting on $\mathcal{C}_c^\infty(X^\circ).$
(This stipulation is important only in dimension two, where $\Lap$ is not
essentially self-adjoint.)

We now can (roughly) state the following:
\begin{theorem}[Melrose--Wunsch \cite{Melrose-Wunsch1}]
Singularities for solutions to the wave equation propagate along diffractive geodesics; strictly
diffractive geodesics generically propagate \emph{weaker}
singularities than geometric geodesics.
\end{theorem}
The genericity condition is that the incident singularities not be
precisely \emph{focused} on the cone tip and applies, e.g., to Cauchy data
that are conormal with respect to a manifold that is at most simply
tangent to the hypersurfaces at constant distance from a cone tip.  In
this case---and in particular for the fundamental solution---we
find that \emph{the diffracted wave for the fundamental solution is
  $(n-1)/2-\epsilon$ derivatives smoother than the main wavefront,}
where $n$ is the dimension of $X.$

We remark that this result has been subsequently generalized to cover
the cases of manifolds with incomplete edge singularities \cite{MVW1},
as well as manifolds with corners \cite{Va:04}, \cite{MeVaWu:13}.

The rest of this paper is essentially applications and extensions of this result in various contexts.

\section{Local energy decay on conic manifolds with Euclidean ends}\label{section:BW}

Consider now a noncompact $n$-manifold $X$ with ends that are
Euclidean.  We will consider solutions to the wave equation
$$
\Box u=0
$$
on $X$ with compactly supported Cauchy data in the energy space.

If $X$ is a smooth manifold, it has long been known that the decay of
local energy can be obstructed by the trapping of geodesics; recall
that a geodesic is said to be forward- or backward-\emph{trapped} if
it remains in a compact set as $t \to \pm \infty.$ Classic work of
Lax--Philips \cite{Lax-Phillips1} and Morawetz \cite{Morawetz:Decay}
shows that, for odd $n,$ absence of trapping implies exponential local
energy decay; on the other hand, results starting with Ralston
\cite{Ralston:Localized} show that trapping of rays implies that
exponential local energy decay cannot hold.  The usual line of
reasoning in obtaining such estimates involves obtained bounds on the
\emph{cutoff resolvent}
$$
\chi (\Lap-\lambda^2)^{-1} \chi,\quad \chi \in \mathcal{C}_c^\infty.
$$
It is well known that in odd dimensions this operator can be
meromorphically continued from $\Im \lambda>0$ to $\CC,$ and its poles
are known as \emph{resonances}.  Exponential local energy decay is
then obtained by showing that no resonances lie in some \emph{strip}
$\Im \lambda>-\nu,$ $\nu>0$ (and that the resolvent has an upper bound
with polynomial growth in this strip).

The situation with conic manifolds is thus interesting for the
following reason: as soon as we have more than one cone point (or,
indeed,\footnote{The author is grateful to Yves Colin de Verdi\`ere
  for pointing out this possibility. In practice, it seems hard to create an
  interesting example of a non-simply connected manifold where the
  \emph{only} trapping is a strictly diffractive geodesic of this form.  On the other hand one may probably add a complex absorbing
  potential to the problem to destroy other trapping and create
  non-simply connected examples.} at least one
cone point if the manifold is non-simply connected) there must be trapped
diffractive geodesics: we can simply continue traversing geodesics
connecting the various cone points.  An example of particular interest
is (the double of) a domain exterior to one or more polygons in $\RR^2$: diffractive
geodesics can move along edges of one polygon and also along lines
connecting vertices of two different polygons.

To what degree, one wonders, does this obstruct energy decay?  The
following theorem (which answers affirmatively a conjecture of
Chandler-Wilde--Graham--Langdon--Spence~\cite{CWGLS:2012} for
polygonal exterior domains) shows that the obstruction is very minor:
\begin{theorem}[Baskin--Wunsch \cite{BaWu:13}]\label{theorem:BaWu}
Assume that no three cone points in $X$ are collinear and no two are
conjugate.  Assume that geodesics missing the cone points escape to
infinity at a uniform rate.

For $\chi \in \CcI(X),$ there exists $\delta>0$ such that the cut-off resolvent
$$
\chi (\Lap-\lambda^2)^{-1}\chi
$$
can be analytically continued from $\Im \lambda >0$ to the region
$$
\Im \lambda >-\rho \log \smallabs{\Re \lambda},\ \smallabs{Re \lambda} >\rho^{-1}
$$
and for some $C,T>0$ enjoys the estimate
$$
\norm{\chi (\Lap-\lambda^2)^{-1}\chi}_{L^{2}\to L^{2}} \leq C
\smallabs{\lambda}^{-1} e^{T\smallabs{\Im \lambda}}
$$
in this region.
\end{theorem}
We contrast this with the the standard result for smooth non-trapping
perturbations of Euclidean space.  In that case the methods of Vainberg
\cite{Vainberg:Asymptotic} and Lax--Phillips \cite{Lax-Phillips1} yield precisely the \emph{same}
resolvent estimate on $\RR$ and a slightly stronger result on
resonance-free regions: \emph{any} region of the form $\Im \lambda>-\rho \smallabs{\log \Re
\lambda}$ is free of resonances outside a large disc.  Thus
the effect of diffractive trapping by cone points is extremely
weak.  Previous results in this direction include energy decay results of \cite{Cheeger-Taylor2}, Section
6, in certain special cases of conic singularities; analogous results
for multiple inverse square potentials were previously proved by
Duyckaerts \cite{Duyckaerts1}.  Burq \cite{Burq:Coin} gave a
precise description of the resonances in the closely related case of
two convex analytic domains in the plane, one of which has a corner
facing the other.  The
diffractive trajectory here bounces back and forth between the corner
and the other obstacle, and Burq showed the associated resonances lie
along a family of logarithmic curves.

We now briefly describe some results on evolution equations that follow from Theorem~\ref{theorem:BaWu}.
We let $\dom_s$ denote the
domain of $\Lap^{s/2}$ (hence locally just $H^s$ away from cone points) and let
$\sin t\sqrt\Lap/\sqrt\Lap$ be the wave propagator.  Let $\chi$ equal $1$
on the set where $X$ is not isometric to $\RR^n.$  In odd dimensions,
the resolvent is a meromorphic function of $\lambda\in \CC$ (with no
difficulties at $\lambda=0$) so in this case
Theorem~\ref{theorem:BaWu} shows that there are only finitely many
resonances in any horizontal strip in $\CC.$  This enables us to show
the following by a contour deformation argument:
\begin{corollary}\label{corollary:resexp}
Let $n$ be odd. Under the assumptions of Theorem~\ref{theorem:BaWu}, for all $A>0,$
small $\ep>0,$ $t>0$ sufficiently large,  and $f \in
\dom_1,$
$$
\chi \frac{\sin t\sqrt{\Lap}}{\sqrt{\Lap}} \chi f=  \sum_{\substack{\lambda_j \in \Res(\Lap) \\ \Im \lambda
  >-A}}\sum_{m=0}^{M_j} e^{-it \lambda_j} t^m w_{j,m} + E_A(t) f
$$
where the sum is of resonances of $\Lap,$ i.e.\ over the poles of the
meromorphic continuation of the resolvent, and the $w_{j,m}$ are the
associated resonant states corresponding to $\lambda_j.$
The error satisfies
$$
\norm{E_A(t)}_{\dom_{1}\to L^{2}} \leq C_\ep e^{-(A-\ep)t}.
$$

In particular, since the resonances have imaginary part bounded above by a
negative constant, $\chi \frac{\sin t\sqrt{\Lap}}{\sqrt{\Lap}} \chi f$ is
exponentially decaying in this case.
\end{corollary}

Another corollary is a local smoothing estimate for the Schr\"odinger
equation.  As it comes from the resolvent estimate on $\RR,$ this is
again lossless as compared to the situation on free $\RR^n$:
\begin{corollary}
  \label{corollary:local-smoothing}
  Suppose $u$ satisfies the Schr{\"o}dinger equation on $X$:
  \begin{align*}
    i^{-1}\pa_t u(t,z) + \Lap  u(t,z) &= 0 \\
    u(0,z) &= u_{0}(z)\in L^{2}(X)
  \end{align*}
Under the assumptions of Theorem~\ref{theorem:BaWu}, for all $\chi \in C^{\infty}_{c}(X)$, $u$ satisfies the local smoothing estimate without loss:
  \begin{equation*}
    \int_{0}^{T}\norm{\chi u(t) }_{\dom_{1/2}}^{2}\, dt \leq
    C_{T} \norm{u_{0}}_{L^2}^{2}.
  \end{equation*}
\end{corollary}

The elements of the proof of Theorem~\ref{theorem:BaWu} are twofold.
The first step is to show that a \emph{very weak Huygens principle}
holds.  We recall that in nontrapping manifolds, a solution to the
wave equation with compactly supported initial data is eventually
\emph{smooth}---this is the usual ``weak Huygens principle.''  Here we show instead that the solution eventually gets
\emph{as smooth as we like}:
\begin{proposition}\label{proposition:Huygens}
Let $\chi \in \CcI(X).$   For any $s \in \RR,$ there exists $T_s\gg 0$
such that whenever $t>T_s,$
$$
\chi U(t) \chi: H^r \to H^{r+s}
$$
for all $r.$
\end{proposition}

The second part of the theorem is a modification of the celebrated
pa\-ra\-met\-rix construction of Vainberg \cite{Vainberg:Asymptotic} (see
also \cite{Tang-Zworski1}).  This argument in its original form builds
a parametrix for the resolvent out of the fundamental solution to the
wave equation, assuming that the latter satisfies the weak Huygens
principle; the new variant, by contrast, makes the weaker assumption
of the output of Proposition~\ref{proposition:Huygens} and produces a very
slightly weaker result (smaller resonance-free region).

Among the further applications of this line of reasoning is the following
theorem on Strichartz estimates for exterior polygonal domains (joint
work with Baskin and Marzuola) \cite{BaskinMarzuolaWunsch:2014}):
for an exterior polygonal domain where the only trapped geodesics are
strictly diffractive (and where no three vertices are collinear) we
find that the same Strichartz estimates for the Schr\"odinger equation
hold as on Euclidean space (locally in time for Neumann conditions,
and globally for Dirichlet).

\section{The wave trace}

If $X$ is a compact Riemannian manifold without boundary let
$$
(\phi_j, \lambda_j^2)
$$
denote the eigenfunctions and eigenvalues of $\Lap.$
One might like to study
the ``inverse spectral problem'' of using the $\lambda_j$ to
characterize $X$ by forming a useful generating function out of the
$\lambda_j.$ An obvious but not directly useful one might be
$$
\sum_j \delta(\lambda-\lambda_j),
$$
but a much more tractable one is the Fourier transform of this
quantity,
$$
\sum_j e^{-it\lambda_j}.
$$
The utility of this generating function stems from its identification
as
$$
\Tr U(t),
$$
where $$U(t)\equiv e^{-it\sqrt{\Lap}}$$ is the ``half-wave'' evolution
operator, mapping functions on $X$ to (certain) solutions to the wave
equation.  If we can say something about the trace of $U(t)$ in terms
of the geometry of $X,$ we can thus hope to learn something about
spectral geometry.

In the setting of smooth
boundaryless manifolds, we have the following classical results on the
wave trace.  Let
$$
\LSpec (X) =\{0\} \bigcup \big\{\pm \text{lengths of periodic
  geodesics on }X\big\}.
$$
\begin{theorem}[Chazarain \cite{Chazarain1}, Duistermaat--Guillemin
  \cite{Duistermaat-Guillemin1}; cf.\ also Colin de Verdi\`ere \cite{Co:73a}, \cite{Co:73b}]\label{theorem:smoothpoisson}
$$\singsupp \Tr U(t) \subset \LSpec (X).$$
\end{theorem}
This allows one to dream of ``hearing'' lengths of closed geodesics,
but does not rule out the possibility that the allowable singularities
do not, in fact, arise.  The presence of honest singularities is,
however, guaranteed by:
\begin{theorem}[Duistermaat--Guillemin \cite{Duistermaat-Guillemin1}]
Let $L$ be the length of an nondegenerate periodic closed geodesic
$\gamma$ on
$L$ that is isolated in the length spectrum.  Then near $t=L$ we have
$$
\Tr U(t) \sim \frac{L_0}{2\pi} i^{\sigma} \abs{I-P}^{-1/2}(t-L)^{-1},
$$
where
\begin{itemize}
\item $L_0$ is the length of the primitive closed geodesic if 
 $\gamma$ is an iterate of a shorter one.
\item $\sigma$ is the Morse index of the variational problem for a
  periodic geodesic, evaluated at $\gamma.$
\item $P$ is the linearized Poincar\'e map, obtained as the
  linearization at $\gamma$ of the first return map to a hypersurface
  of the phase space, transverse to $\gamma.$
\end{itemize}
\end{theorem}
Note that the nondegeneracy condition in the hypotheses is simply the condition that $I-P$ be nonsingular.

The generalization of Theorem~\ref{theorem:smoothpoisson} to compact
conic manifolds is straightforward: let $$\DLSpec(X)=
\{0\} \bigcup \big\{\pm \text{lengths of periodic
 diffractive geodesics on }X\big\}.$$
\begin{theorem}[Wunsch \cite{Wunsch2}]
On a conic manifold $X,$
$$
\singsupp \Tr U(t) \subset \DLSpec{X}.
$$
\end{theorem}
The singularities at lengths of geodesics in $X^\circ$ are easily seen
to be described by the same formula given by Duistermaat--Guillemin,
but the geodesics interacting through conic points are not so simple.
We consider $\gamma$ a closed, \emph{strictly diffractive} geodesic
undergoing $k$ diffractions and traversing geodesic segments
$\gamma_1,\dots, \gamma_k$ connecting cone points $Y_{i_1}, \dots
Y_{i_k}.$ Recall that the hypothesis that the geodesic be strictly
diffractive means that it interacts with each cone point by entering
and leaving on a pair of geodesics that cannot be uniformly locally
approximated by geodesics in $X^\circ.$ This is generically the case
for all closed geodesics.  Assume further that the length $L$ of
$\gamma$ is isolated in the length spectrum, and make the additional
nondegeneracy hypothesis that no two cone points along the geodesic
are conjugate to one another.
Note that the following was previously known by work of Hillairet
\cite{Hillairet:2005} in the important special case of flat surfaces
with conic singularities (hence in particular for doubles of
polygons).
\begin{theorem}[Ford--Wunsch \cite{1411.6913}]\label{theorem:FoWu}
Near $t=L,$
$$
\Tr U(t) \sim \int e^{i(t-L)\xi} a(\xi) \, d\xi
$$
where
\begin{equation}\label{symbol}
  a(\xi) \sim L_0 \cdot (2\pi)^{\frac{kn}{2}} \,  e^{\frac{ik(n-3)\pi}{4}} \,
  \chi(\xi) \, \xi^{-\frac{k(n-1)}{2}} \prod_{j=1}^k i^{-m_{\gamma_j}} \, \D_j
  \, \W_j \ \text{as $|\xi| \to \infty$}. 
\end{equation}
\end{theorem}
Here, 
$\chi\in \CI(\RR)$ is $1$ for $\xi>1$ and $0$ for
  $\xi<0.$  
Note that the power of $\xi$ is such that we obtain greater smoothness
as the number of diffractions increases.  The leading order
singularity as a function of $t$ is proportional to $(t-L+i0)^{-1+k(n-1)/2}$
  (but is multiplied by $\log(t-L+i0)$ if the power is an integer).

As before $L_0$ denotes the length of the
``primitive'' geodesic if $\gamma$ is an iterate of a shorter one.
The integers $m_{\gamma_j}$ are simply the Morse indices of the
variational problems associated to traveling from one cone point to
the next, evaluated at $\gamma_j.$

We will now explain the factors $\D_j$ and $\W_j.$

The terms $\D_j$ are associated to the diffractions through each
successive cone point $Y_{i_j}.$  They are constructed as follows.
Each cone point $Y_{i_j}$ is equipped with a metric $h_{i_j}\equiv
h\restrictedto_{Y_{i_j}}.$  It thus has a Laplace-Beltrami operator
$\Lap_{i_j}$ and we may use the functional calculus to take functions
of this operator.  In particular, let $$\nu_{i_j} \equiv \sqrt{
  \Delta_{Y_{i_j}} + \left( \frac{2-n}{2} \right)^2 }.$$
We then form the operator family
$$
e^{-it\nu_{i_j}}: L^2(Y_{i_j})\to L^2(Y_{i_j}).
$$
This is essentially a ``half Klein Gordon propagator'' on the link of
the cone point (i.e., a boundary component).  Now let $\kappa(\bullet)$ denote the Schwartz
kernel of an operator.  Supposing that the diffractive geodesic
$\gamma$ enters $Y_{i_j}$ at the point $y$ and leaves from point $y',$
we set
$$
\D_j \equiv \kappa(e^{-i\pi \nu_{i_j}})[y,y'].
$$
The propagator kernel is of course not continuous in general, however note
that the strictly diffractive nature of the geodesic ensures that $y$
and $y'$ are not connected by a geodesic of length $\pi$ in the link,
which in turn precisely ensures, by propagation of singularities, that
the Schwartz kernel of the time-$\pi$ Klein Gordon propagator is
smooth near $(y,y'),$ hence the evaluation of this distribution makes
sense.

Now we turn to $\W_j$.  These quantities are associated to the
geodesic segments $\gamma_j$ connecting successive cone points.  They
are best described in terms of Jacobi fields, but can also be
viewed as a proxy for a quantity involving the derivative of the
expenential map, hence a substitute for the term involving the
Poincar\'e map in the Duistermaat--Guillemin formula.   Note that the
exponential map from one cone point to the next does not make sense,
since any small perturbation of the geodesic $\gamma_j$ will miss the
next cone point entirely rather than simply hitting it at a different point.
Correspondingly, if we let $\bJ$ be a set of Jacobi fields that are
orthonormal to $\gamma_j$ and at $\gamma_j(0)$ give an orthonormal
basis of $TY_{i_j}$ then $\bJ$ becomes \emph{singular} as we approach
the end of $\gamma_j$ at $Y_{i_{j+1}}.$  On the other hand, the metric
is also singular at cone points, in the sense that it vanishes on
$TY,$ so we can nonetheless make sense of the determinant
$$
\det_g \bJ\restrictedto_{Y_{i_{j+1}}}.
$$
Then we have
$$
\W_j \equiv\big\lvert\det_g \bJ\restrictedto_{Y_{i_{j+1}}}\big\rvert^{-1/2}.
$$
This quantity can be made to look more like the derivative of an
exponential map as follows: we set
\begin{equation}\label{thetaj}
\Theta_j=(\length(\gamma_j)^{-(n-1)})\big\lvert\det_g \bJ\restrictedto_{Y_{i_{j+1}}}\big\rvert.
\end{equation}
Consider the case in which $Y_j$ is a ``fictitious'' cone point
obtained by blowing up a smooth point $p_0$ on a manifold.  Then Jacobi vector
fields tangent to $Y_j$ are obtained as lifts under the blow-down map
of Jacobi fields vanishing at $p_0,$ and 
$\Theta_j$ becomes
a standard expression
for $\det D\exp_{Y_j} (\bullet)$ in terms of Jacobi fields, at least
when evaluated in $X^\circ$ (cf.\
\cite{Be:77}): in that case we simply have
$$
\Theta_j=\det_g \big\lvert D \exp_{Y_j}(\bullet)\big\rvert.
$$
Since $\W_j=\length(\gamma_j)^{-(n-1)/2} \Theta_j^{-1/2}$ we
recover the relationship with the exponential map in the case of a
trivial cone point.

In rough outline, the proof of Theorem~\ref{theorem:FoWu} goes as
follows.  We know explicitly what the wave propagator look like on a
model \emph{product cone} $\RR_+\times Y_j$ endowed with the scale
invariant metric $dx^2 +x^2 h_0(y,dy)$---this is a computation of
Cheeger--Taylor \cite{Cheeger-Taylor1}, \cite{Cheeger-Taylor2}
involving bravura use of the Hankel transform.  In particular, we
can evaluate the symbol of the diffracted wavefront explicitly in that
case.  More generally, in \cite{Melrose-Wunsch1} the author and
Melrose prove that near a cone point, the diffracted front of the wave
propagator is guaranteed to be a conormal distribution.  The first new
step is therefore to show that in the non-product case, the principal
symbol of the diffracted front is still, modulo adjustments involving
comparing half-densities on the two spaces, given by the same
expression as in the product case where we use the model metric
$dx^2+x^2 h\restrictedto_{x=0} (y,dy).$ This involves comparing the two
propagators and showing that the difference
between model and exact propagators can be estimated by a
\emph{Morawetz inequality} near the cone tip.

Having understood the effect of a single diffraction, we then proceed
as follows.  We take a microlocal partition of unity $A_j$ on $X,$
where for technical reasons the $A_j$ are restricted to be simply
cutoff functions near each boundary component $Y_i$ but are otherwise
fully localized in phase space.  
We then decompose the wave trace as follows: fix
small times $t_j$ with $\sum t_j=T.$  Then by cyclicity of the trace
$$
\Tr U(t) = \sum_{i_0,\dots, i_{N}}\Tr \sqrt{A_{i_0}} U(t-T)A_{i_1} U(t_1) A_{i_2}\dots
A_{t_{N}} U(t_N) \sqrt{A_{i_0}}.
$$
By propagation of singularities, this term is guaranteed to be trivial
unless there is a diffractive geodesic successively passing through
the microsupports of the $A_{i_j}$'s, hence we may throw away most of
this sum.  The remaining terms are then computed by a stationary
phase computation, gluing together the propagators for ``free''
propagation through $X^\circ$ with those for the diffractive
interaction with cone points (this was the same strategy previously
used by Hillairet in \cite{Hillairet:2005} as well as by the author in
\cite{Wunsch2}).

\section{Lower bounds for resonances}

While $\Tr U(t)$ only makes sense (even distributionally) on a \emph{compact}
manifold, if we return to the setting of Section~\ref{section:BW}
where we have a \emph{noncompact} manifold with Euclidean ends, we may still
make sense of an appropriately \emph{renormalized} wave trace, and use
the diffractive trace formula (Theorem~\ref{theorem:FoWu}) to obtain
lower bounds on resonances.

In odd dimensions, we let $\wavegen$ denote the generator of
the wave group, and hence $e^{t \wavegen}$ the wave group itself;
likewise we let $\wavegen_0$ be the generator of the wave group on Euclidean space.
 We
then have the trace formula
\begin{equation}\label{traceres}
\Tr (e^{t \wavegen} -e^{t \wavegen_0}) = \sum_{\lambda_j \in \Res} e^{-i \lambda_j t},\ t
> 0
\end{equation}
where the sum is over the resonances, counted with multiplicity (see
e.g. \cite{Sjostrand-Zworski5} for the details of how to makes sense of
this difference of operators in a wide variety of contexts).
This result in various settings was first proved by Bardos-Guillot-Ralston \cite{Bardos-Guillot-Ralston1}, Melrose
\cite{MR83j:35128}, and Sj\"ostrand-Zworski \cite{Sjostrand-Zworski5};
an analogous result in even dimensions can be found in \cite{Zw:99}.

Now if we can actually guarantee the existence of singularities in the (renormalized) wave
trace, a Tauberian theorem of Sj\"ostrand-Zworski
\cite{Sjostrand-Zworski3} allows us to deduce from \eqref{traceres} in
a lower bound on the number of resonances in logarithmic
regions in $\CC.$  Fortunately, Theorem~\ref{theorem:FoWu} applies
equally well in this context, and we obtain a lower bound on the
number of resonances as follows.
Let
$$
N_\rho(r) =\# \{\text{Resonances in } \smallabs{\lambda}<r,\
\Im \lambda \geq -\rho \log \smallabs{\Re{\lambda}}\}.
$$
Then we have:
\begin{theorem}[Hillairet--Wunsch]\label{theorem:HiWu}
Under the geometric assumptions of Theorem~\ref{theorem:BaWu}, let $L$
be the length of a closed, strictly diffractive geodesic $\gamma$
undergoing $k$ diffractions.  Assume, in
the notation of Theorem~\ref{theorem:FoWu}, that all the diffraction
coefficients $\D_j$ are nonzero along $\gamma;$ assume also that there are no closed diffractive geodesics
beside iterates of this one having length in $L\NN.$ Then for all $\ep>0,$
$$
N_\rho(\rho)\geq C_{\rho,\epsilon} r^{1-\epsilon}
$$
provided
$$
\rho > \frac{(n-1)k}{2L}
$$
\end{theorem}
A detailed proof, which simply consists of using the trace formula
(Theorem~\ref{theorem:FoWu}) in \eqref{traceres} together with the
Tauberian theorem of \cite{Sjostrand-Zworski3}, can be found in
\cite{Ga:15}.  Note that the bound on $\rho$ written here is that
which we obtain by considering the whole sequence of singularities of
the wave trace obtained by considering arbitrary \emph{iterates} of
the geodesic $\gamma.$ We remark that the distinction between the
trace of the full wave group and $\Tr U(t)$ is immaterial for this
purpose since the former is twice the real part of the latter, and it
is not difficult to verify from examination of \eqref{symbol} that the
singularities arising from iterates of a given geodesic cannot all be
purely imaginary.

The optimal $\rho$ here is generally obtained by choosing $\gamma$
to be the geodesic that
traverses the longest geodesic segment connecting a pair of distinct cone
points, back and forth (assuming the diffraction coefficients are
nonvanishing).  If
$D_{\text{max}}$ denotes the greatest distance between a pair of cone points, then we
have a closed geodesic of length $2 D_{\text{max}}$ with $k=2,$ and we
obtain the bound
$$
\rho > \frac{(n-1)}{2D_{\text{max}}}.
$$

Remarkably, this theorem is essentially sharp, as was shown by
Galkowski, who has produced an effective version of the Vainberg
argument previously employed in \cite{BaWu:13}:
\begin{theorem}[Galkowski \cite{Ga:15}]
Let $D_{\text{max}}$ be the greatest distance between two cone points.  For any
$\epsilon>0$ the
constant $\rho$ in Theorem~\ref{theorem:BaWu} can be taken to be
$(n-1)/(2D_{\text{max}})-\epsilon,$ i.e.\ $N_\rho(r)$ is \emph{bounded} for all $\rho<(n-1)/2D_{\text{max}}.$
\end{theorem}
Since $N_\rho$ is bounded for $\rho<(n-1)/2D_{\text{max}}$ and (subject to the
nondegeneracy hypotheses of Theorem~\ref{theorem:HiWu}) almost linearly
growing for $\rho>(n-1)/2D_{\text{max}},$ we find that in any set near the
critical curve $\Im\lambda=-((n-1)/2D_{\text{max}}) \log \smallabs{\Re
  \lambda}$ of the form
$$
\big(-\frac{n-1}{2D_{\text{max}}}-\ep\big) \log \smallabs{\Re
  \lambda}<\Im \lambda<\big(-\frac{n-1}{2D_{\text{max}}}+\ep\big)
\log \smallabs{\Re \lambda},\quad \abs{\lambda}>\ep^{-1}
$$
there are infinitely many resonances.  The intuition behind the
importance of the longest geodesic connecting two cone points is that
repeatedly traversing this segment back and forth is the way in which
a trapped singularity can diffract \emph{least frequently}.  Since each
diffraction loses considerable energy owing to the smoothing effect of
diffraction, a resonant state propagating back and forth along this
geodesic is the one that loses energy to infinity at the
slowest rate.

\bibliographystyle{plain}
\bibliography{all}
\end{document}